\documentclass{amsart}
\usepackage{amsfonts,amsmath,amsthm,amssymb}
\newtheorem{theorem}{Theorem}

\newtheorem{corollary}{Corollary}
\newtheorem{proposition}{Proposition}
\title{Torsion units in integral group ring of Higman-Sims simple group}
\author{V.A.~Bovdi,  A.B.~Konovalov}
\date{}

\address{
V.A.~Bovdi\newline
Institute of Mathematics, University of
Debrecen\newline P.O.  Box 12, H-4010 Debrecen,
Hungary\newline Institute of Mathematics and
Informatics, College of Ny\'\i regyh\'aza\newline
S\'ost\'oi \'ut 31/b, H-4410 Ny\'\i regyh\'aza, Hungary}
\email{vbovdi@math.klte.hu}

\address{
A.B.~Konovalov
\newline School of Computer Science, University of St Andrews,
\newline Jack Cole Building, North Haugh, St Andrews, Fife, KY16 9SX, Scotland}
\email{konovalov@member.ams.org}


\thanks{
The research was supported by OTKA grants No.T 037202, No.T 038059}

\dedicatory{Dedicated to  Professor  B\'ela Cs\'ak\'any  on his
75th birthday}


\begin{document}

\begin{abstract}
Using the Luthar--Passi method, we investigate the classical
Zassenhaus conjecture for the normalized unit group of the
integral group ring of the Higman-Sims simple sporadic group
$\verb"HS"$. As a consequence, we confirm the Kimmerle's
conjecture on prime graphs for this sporadic group.
\end{abstract}

\maketitle

\section{Introduction  and main results}
\label{Intro}

Let $V(\mathbb Z G)$ be  the normalized unit group of the integral
group ring $\mathbb Z G$ of a finite group $G$. One of most
interesting conjectures in the theory of integral group ring is
the conjecture {\bf (ZC)} of  H.~Zassenhaus \cite{Zassenhaus},
saying that every torsion unit $u\in V(\mathbb ZG)$ is conjugate
to an element in $G$ within the rational group algebra $\mathbb Q
G$.

For finite simple groups, the main tool of  the investigation of
the Zassenhaus conjecture is the Luthar--Passi method, introduced
in \cite{Luthar-Passi} to solve the (ZC) for $A_{5}$. Later in
\cite{Hertweck1} M.~Hertweck extended  and applied it for the
investigation of the Zassenhaus conjecture for $PSL(2,p^n)$. The
method proved to be useful for groups containing non-trivial normal
subgroups as well. We refer to
\cite{Bovdi-Hofert-Kimmerle,Bovdi-Konovalov,Hertweck2, Hertweck1,
Hertweck3, Hofert-Kimmerle} for  recent results. Related results
can be found in \cite{Artamonov-Bovdi,Luthar-Trama} and
\cite{Bleher-Kimmerle,Kimmerle}. In the latter papers weakened
versions  of the (ZC) were conjectured.

In order to state one of these we introduce some notation. By $\#
(G)$ we denote the set of all primes dividing the order of $G$.
The Gruenberg--Kegel graph (or the prime graph) of $G$ is the graph
$\pi (G)$ with vertices labelled by the primes in $\# (G)$ and
there is an edge from $p$ to $q$ if and only if there is an element of order
$pq$ in the group $G$. In \cite{Kimmerle} W.~Kimmerle  proposed
the following:

\centerline {{\bf Conjecture (KC)}: if $G$ is a
finite group then $\pi (G) =\pi (V(\mathbb Z G))$.}

Obviously,  the Zassenhaus conjecture {\bf (ZC)} implies the
Kimmerle conjecture {\bf (KC)}. In \cite{Kimmerle} it was shown, that
{\bf (KC)} holds for finite Frobenius and solvable groups. We remark
that with respect to the so-called $p$-version of the Zassenhaus
conjecture the investigation of Frobenius groups was completed by
M.~Hertweck and the first author in \cite{Bovdi-Hertweck}. In
\cite{Bovdi-Konovalov, Bovdi-Konovalov-M23,
Bovdi-Konovalov-Linton, Bovdi-Konovalov-Siciliano}, {\bf (KC)} was also confirmed for
certain Mathieu sporadic simple groups, and in \cite{Bovdi-Jespers-Konovalov}
-- for some Janko sporadic simple groups.

In this paper we continue these investigations for the Higman-Sims
simple sporadic group $\verb"HS"$. The main result provides
information about the possible torsion units in $V(\mathbb Z
\verb"HS")$. An immediate consequence is a positive answer to {\bf
(KC)} for $\verb"HS"$.

In order to state the result we need to introduce some notation.
Let $G$ be a group. Let $\mathcal{C} =\{ C_{1}, \ldots, C_{nt},
\ldots \}$ be  the collection of all conjugacy classes of $G$,
where the first index denotes the order of the elements of this
conjugacy class and $C_{1}=\{ 1\}$. Suppose $u=\sum \alpha_g g \in
V(\mathbb Z G)$ has finite order $k$. Denote by
$\nu_{nt}=\nu_{nt}(u)=\varepsilon_{C_{nt}}(u)=\sum_{g\in C_{nt}}
\alpha_{g}$, the partial augmentation of $u$ with respect to
$C_{nt}$. From the Berman--Higman Theorem 
(see \cite{Berman} and \cite{Sandling}, Ch.5, p.102) 
one knows that
$\nu_1 =\alpha_{1}=0$ and
\begin{equation}\label{E:1}
\sum_{C_{nt}\in \mathcal{C}} \nu_{nt}=1.
\end{equation}
Hence, for any character
$\chi$ of $G$, we get that $\chi(u)=\sum
\nu_{nt}\chi(h_{nt})$, where $h_{nt}$ is a
representative of a conjugacy class $ C_{nt}$.

The main result is the following.

%
%

\begin{theorem}\label{T:1}
Let $G$ denote the Higman-Sims simple sporadic  group $\verb"HS"$.
Let $u$ be a torsion unit of $V(\mathbb ZG)$ of order $|u|$.
Denote by $\frak{P}(u)$ the tuple
\[
\begin{split}
(\nu_{2a}, \nu_{2b}, \nu_{3a}, \nu_{4a},& \nu_{4b}, \nu_{4c}, \nu_{5a},
\nu_{5b}, \nu_{5c}, \nu_{6a}, \nu_{6b}, \nu_{7a},\\
&\nu_{8a}, \nu_{8b}, \nu_{8c}, \nu_{10a}, \nu_{10b}, \nu_{11a},
\nu_{11b}, \nu_{12a}, \nu_{15a}, \nu_{20a}, \nu_{20b}) \in \mathbb
Z^{23}
\end{split}
\]
of partial augmentations of $u$ in $V(\mathbb ZG)$. The following
properties hold.

\begin{itemize}

\item[(i)] 
There is no elements of orders $14$, $21$, $22$, $33$, $35$, $55$ and 
$77$ in $V(\mathbb ZG)$.  Equivalently, 
if  $|u|\not\in\{24, 30, 40, 60, 120\}$, then $|u|$ coincides
with the order of some  $g\in G$.

\item[(ii)]If $|u| \in \{3,7 \}$, then $u$ is rationally conjugate
to some $g\in G$.

\item[(iii)] If $|u|=2$, the tuple of the  partial augmentations
of $u$ belongs to the set
\[
\begin{split}
\big\{\; \frak{P}(u) \mid (\nu_{2a},\nu_{2b}) \in \{ \; (0,
1), (-2, 3), (2, -1), (1, 0), (3, -2), (-1, 2) \; \},&\\
\nu_{kx}=0,\; kx\not\in\{2a,2b\}\; &\big\}.
\end{split}
\]

\item[(iv)]
 If $|u|=5$, the tuple of the partial augmentations
of $u$ belongs to the set
\[
\begin{split}
\big\{\; \frak{P}(u) \mid (\nu_{5a},\nu_{5b},\nu_{5c}) \in \{ \; (
-2, -1, 4 ),\; ( -1, -1, 3 ),\; ( 0, -1, 2 ),\; ( 1, -1, 1 ),&\\
( 1, 3, -3 ),\; ( 0, 3, -2 ),\; ( -3, 0, 4 ),\; ( -2, 0, 3 ),\; (
1,0, 0 ),\;( 1, 4, -4 ), &\\
( -1, 0, 2 ),\; ( 0, 0, 1 ),\; ( 0, 2, -1), \; (-1, 2, 0 ),\; ( 1,
2, -2 ),\; ( 0, 1, 0 ),&\\
( -2, -2, 5 ),\; ( -1, -2, 4 ),\; ( 0, -2, 3 ),\; ( 1, -2, 2 ),\;
( 1, 1, -1 ),&\\
\; ( -2, 1, 2 ),\; ( -1, 1, 1 ) \; \},\qquad  \nu_{kx}=0,\;
kx\not\in\{5a,5b,5c\} & \; \big\}.
\end{split}
\]

\item[(v)] If $|u|=11$, the tuple of the  partial augmentations of
$u$ belongs to the set
\[
\begin{split}
\big\{\; \frak{P}(u) \mid (\nu_{11a},\nu_{11b})) \in \{ \; (5, -4
),\; ( 4, -3 ),\; ( -2, 3 ),\; ( 2, -1 ),\; ( -3, 4 ),&\\
( -4, 5
),\; ( 1, 0 ),\; ( 3, -2 ),\; ( -1, 2 ),\; ( 0, 1 )\; \}, & \\
\nu_{kx}=0,\quad  kx\not\in\{11a,11b\} & \; \big\}.
\end{split}
\]

\end{itemize}
\end{theorem}

%
%

\begin{corollary} If $G$ is the Higman-Sims sporadic group, then
$\pi(G)=\pi(V(\mathbb ZG))$.
\end{corollary}

\section{Preliminaries}

The following result relates the solution of
the Zassenhaus conjecture to vanishing of 
partial augmentations of torsion units.

\begin{proposition}\label{P:5}
(see \cite{Luthar-Passi} and
Theorem 2.5 in \cite{Marciniak-Ritter-Sehgal-Weiss})
Let $u\in V(\mathbb Z G)$
be of order $k$. Then $u$ is conjugate in $\mathbb
QG$ to an element $g \in G$ if and only if for
each $d$ dividing $k$ there is precisely one
conjugacy class $C$ with partial augmentation
$\varepsilon_{C}(u^d) \neq 0 $.
\end{proposition}

The next result yield that several partial augmentations are zero.

\begin{proposition}\label{P:4}
(see \cite{Hertweck2}, Proposition 3.1;
\cite{Hertweck1}, Proposition 2.2)
Let $G$ be a finite
group and let $u$ be a torsion unit in $V(\mathbb
ZG)$. If $x$ is an element of $G$ whose $p$-part,
for some prime $p$, has order strictly greater
than the order of the $p$-part of $u$, then
$\varepsilon_x(u)=0$.
\end{proposition}

Another important restriction on the partial augmentations is
given by the following result.

\begin{proposition}\label{P:1}
(see \cite{Luthar-Passi,Hertweck1}) Let either $p=0$ or $p$ is a
prime divisor of $|G|$. Suppose that $u\in V( \mathbb Z G) $ has
finite order $k$ and assume that $k$ and $p$ are coprime when 
$p\neq 0$. If $z$ is a complex primitive $k$-th root of unity and $\chi$
is either a classical character or a $p$-Brauer character of $G$
then, for every integer $l$, the number
$$
\mu_l(u,\chi, p )=\textstyle\frac{1}{k} \sum_{d|k}Tr_{ \mathbb Q
(z^d)/ \mathbb Q } \{\chi(u^d)z^{-dl}\}
$$
is a non-negative integer.
\end{proposition}

Note that if $p=0$, we will use the notation $\mu_l(u, \chi, * )$
for $\mu_l(u,\chi, 0)$.

Finally, we shall use the well-known bound for
orders of torsion units.

\begin{proposition}\label{P:2}  (see \cite{Cohn-Livingstone})
The order of a torsion element $u\in V(\mathbb ZG)$
is a divisor of the exponent of $G$.
\end{proposition}

\section{Proof of the Theorem}

In this section we denote by $G$ the Higman-Sims simple sporadic
group $\verb"HS"$. It is well known \cite{Atlas,GAP} that
$$
|G|=44352000 = 2^{9} \cdot 3^{2} \cdot 5^{3} \cdot 7 \cdot 11
\quad \text{and}\quad exp(G)=9240 = 2^{3} \cdot 3 \cdot 5 \cdot 7
\cdot 11.
$$
The character table of $G$, as well as the Brauer character tables
(denoted by $\mathfrak{BCT}{(p)}$, where $p \in \{2, 3, 5, 7, 11\}$) 
can be found by the computational algebra system GAP \cite{GAP}, 
which derives its data from \cite{Atlas,AtlasBrauer}. 
Throughout the paper we will use the notation of GAP Character
Table Library for the characters and conjugacy classes of the 
Higman-Sims group $\verb"HS"$.

From the structure of the group $\verb"HS"$ we know that it
possesses elements of orders $2$, $3$, $4$, $5$, $6$, $7$, $8$,
$10$, $11$, $12$, $15$ and $20$. We begin our investigation with
units of orders $2$, $3$, $5$, $7$ and $11$. We do no treat the
remaining cases ($4$, $6$, $8$, $12$, $15$, $20$), because in
these cases the computation is quite complex. Since by Proposition
\ref{P:2}, the order of each torsion unit divides the exponent of
$G$, it remains to consider in addition only the units of orders 
$14$, $21$, $22$, $24$, $30$, $33$, $35$, $40$, $55$, $60$, $77$ and $120$. 
Now we omit five remaining cases: $24$, $30$, $40$, $60$ and $120$, 
since these cases are computationally too complicated. We can prove
that the order of an unit can not be equal to $14$, $21$, $22$, 
$33$, $35$, $55$ or $77$. 

Thus, in this paper we are going to treat 
the cases when the order of $u$ is: $2$, $3$, $5$, $7$, $11$, $14$, 
$21$, $22$, $33$, $35$, $55$ and $77$.

\noindent $\bullet$ Let $|u| \in \{3,7\}$. Since there is only one
conjugacy class in $G$ consisting of elements or order $|u|$, this
case follows at once from Proposition \ref{P:4}. Thus, for units
of orders $3$ and $7$ we obtained that there is precisely one
conjugacy class with non-zero partial augmentation. Proposition
\ref{P:5} then yields part (ii) of the Theorem.

\noindent$\bullet$ Let $u$ be an involution. By (\ref{E:1}) and
Proposition \ref{P:4} we get $\nu_{2a}+\nu_{2b}=1$.  Applying
Proposition \ref{P:1} to the character $\chi_{2}$ we get the
following system of inequalities
\[
\begin{split}
\mu_{0}(u,\chi_{2},*) & = \textstyle \frac{1}{2} (6 \nu_{2a} - 2 \nu_{2b} + 22) \geq 0; \\ 
\mu_{1}(u,\chi_{2},*) & = \textstyle \frac{1}{2} (-6 \nu_{2a} + 2 \nu_{2b} + 22) \geq 0. \\ 
\end{split}
\]
From the requirement that all $\mu_i(u,\chi_{j},p)$ must be
non-negative integers it can be deduced that $(\nu_{2a},\nu_{2b})$
satisfies the conditions of  part (iii) of the Theorem.

\noindent $\bullet$ Let $u$ be a unit of order $5$. By (\ref{E:1})
and Proposition \ref{P:4} we obtain that
$$\
\nu_{5a}+\nu_{5b}+\nu_{5c}=1.
$$
Put $t_1 = 3 \nu_{5a} - 2 \nu_{5b} - 2 \nu_{5c}$ and $t_2 = 2
\nu_{5a} - 3 \nu_{5b} + 2 \nu_{5c}$.  Again applying Proposition
\ref{P:1} to the characters $\chi_2$, $\chi_3$ and  $\chi_4$,   we
obtain  the system of inequalities
\[
\begin{split}
\mu_{0}(u,\chi_{2},*) & = \textstyle \frac{1}{5} (-4t_1 + 22) \geq
0; \qquad \mu_{1}(u,\chi_{2},*)  = \textstyle \frac{1}{5} (t_1 + 22) \geq 0; \\ 
\mu_{0}(u,\chi_{3},*) & = \textstyle \frac{1}{5} (4t_2 + 77) \geq
0; \qquad \mu_{1}(u,\chi_{3},*) = \textstyle \frac{1}{5} (-t_2 + 77) \geq 0; \\ 
&\mu_{0}(u,\chi_{3},2)  = \textstyle \frac{1}{5} (24 \nu_{5a} - 16 \nu_{5b} + 4 \nu_{5c} + 56) \geq 0; \\ 
&\mu_{0}(u,\chi_{4},2)  = \textstyle \frac{1}{5} (28 \nu_{5a} + 8 \nu_{5b} - 12 \nu_{5c} + 132) \geq 0; \\ 
&\mu_{0}(u,\chi_{3},3)  = \textstyle \frac{1}{5} (-4 \nu_{5a} + 16 \nu_{5b} - 4 \nu_{5c} + 49) \geq 0, \\ 
\end{split}
\]
that has only $23$ integer solutions
$(\nu_{5a},\nu_{5b},\nu_{5c})$ (they are listed in the part (iv)
of the Theorem) such that all $\mu_i(u,\chi_{j},p)$ are
non-negative integers.

\noindent$\bullet$ Let $u$ be a unit of order $11$. By (\ref{E:1})
and Proposition \ref{P:4} we have
$$
\nu_{2a}+\nu_{3a}+\nu_{6a}=1.
$$
Applying Proposition \ref{P:1} to the characters in
$\mathfrak{BCT}{(3)}$ and $\mathfrak{BCT}{(5)}$ we get the system
\[
\begin{split}
\mu_{1}(u,\chi_{14},*) & = \textstyle \frac{1}{11} (6 \nu_{11a} - 5 \nu_{11b} + 896) \geq 0; \\ 
\mu_{2}(u,\chi_{14},*) & = \textstyle \frac{1}{11} (-5 \nu_{11a} + 6 \nu_{11b} + 896) \geq 0; \\ 
\mu_{1}(u,\chi_{3},3)  & = \textstyle \frac{1}{11} (6 \nu_{11a} - 5 \nu_{11b} + 49) \geq 0; \\ 
\mu_{2}(u,\chi_{3},3)  & = \textstyle \frac{1}{11} (-5 \nu_{11a} + 6 \nu_{11b} + 49) \geq 0; \\ 
\mu_{1}(u,\chi_{9},5)  & = \textstyle \frac{1}{11} (6 \nu_{11a} - 5 \nu_{11b} + 280) \geq 0; \\ 
\mu_{2}(u,\chi_{9},5)  & = \textstyle \frac{1}{11} (-5 \nu_{11a} + 6 \nu_{11b} + 280) \geq 0, \\ 
\end{split}
\]
that has only ten  integer solutions
$(\nu_{11a},\nu_{11b},\nu_{6a})$ listed in  part (v) of the
Theorem.

\noindent$\bullet$ Let $u$ be a unit of order $14$. By (\ref{E:1})
and Proposition \ref{P:4} we have
\begin{equation}\label{E:2}
\nu_{2a}+\nu_{2b}+\nu_{7a}=1.
\end{equation}
Put
\begin{equation}\label{E:3}
{ (\alpha,\beta,\gamma,\delta)=\begin{cases}
(34,22,90,64),\quad &\text{ if }\quad  \chi(u^{7})=\chi(2a);\\
(26,30,78,76),\quad &\text{ if }\quad  \chi(u^{7})=\chi(2b);\\
(10,3,54,100),\quad &\text{ if }\quad  \chi(u^{7})=-2\chi(2a)+3\chi(2b);\\
(42,14,102,52),\quad &\text{ if }\quad  \chi(u^{7})=2\chi(2a)-\chi(2b);\\
(-1,6,114,40),\quad &\text{ if }\quad  \chi(u^{7})=3\chi(2a)-2\chi(2b);\\
(18,38,66,88),\quad &\text{ if }\quad  \chi(u^{7})=-\chi(2a)+2\chi(2b).\\
\end{cases}}
\end{equation}
Additionally, we set
\begin{equation}\label{E:4}
t_1 = 6 \nu_{2a} - 2 \nu_{2b} + \nu_{7a}\qquad  \text{and}\qquad
t_2 = 13 \nu_{2a} + \nu_{2b}.
\end{equation}
Since  $|u^7|=2$, for any character $\chi$ of $G$ we need to
consider six cases, defined by part (iii) of the Theorem. Using
Proposition \ref{P:1}, it is easy to check, that in all six cases
we have the following system of inequalities:
\begin{equation}\label{E:5}
\begin{split}
\mu_{0}(u,\chi_{3},*) & = \textstyle \frac{1}{14} (6 t_2 +\gamma )
\geq 0; \qquad \mu_{7}(u,\chi_{3},*) = \textstyle \frac{1}{14} (- 6 t_2 +\delta) \geq 0. \\ 
\end{split}
\end{equation}
Furthermore,  if\quad  $\chi(u^{7})\in\{\;  \chi(2a),\;
\chi(2b),\; 2\chi(2a)-\chi(2b),\; -\chi(2a)+2\chi(2b)\; \}$,\quad
then again by Proposition \ref{P:1} we get the system:
\begin{equation}\label{E:6}
\begin{split}
\mu_{0}(u,\chi_{2},*) & = \textstyle \frac{1}{14} (6 t_1 + \alpha)
\geq 0;\qquad \mu_{7}(u,\chi_{2},*)  = \textstyle \frac{1}{14} (-6t_1 + \beta) \geq 0. \\ 
\end{split}
\end{equation}
If\quad  $\chi(u^{7})=-2\chi(2a)+3\chi(2b)$\quad  then similarly
as before,  we obtain that
\begin{equation}\label{E:7}
\begin{split}
\mu_{0}(u,\chi_{2},*) & = \textstyle \frac{1}{14} (6 t_1 + \alpha)
\geq 0;\qquad \mu_{2}(u,\chi_{2},*)  = \textstyle \frac{1}{14} (- t_1 + \beta) \geq 0. \\ 
\end{split}
\end{equation}
Finally, if\quad  $\chi(u^{7})=3\chi(2a)-2\chi(2b)$\quad  then by
Proposition \ref{P:1} we have
\begin{equation}\label{E:8}
\begin{split}
\mu_{1}(u,\chi_{2},*) & = \textstyle \frac{1}{14} (t_1 + \alpha)
\geq 0;\qquad \mu_{7}(u,\chi_{2},*)  = \textstyle \frac{1}{14} (-6t_1 + \beta) \geq 0. \\ 
\end{split}
\end{equation}
If we substitute the possible values of
$(\alpha,\beta,\gamma,\delta)$  from (\ref{E:3}) into
(\ref{E:5})-- (\ref{E:8}), then we can compute the possible values
of $t_1$ and $t_2$ in all of the tree cases. Now we substitute back
these values of $t_1$ and $t_2$ into (\ref{E:4}). Then (\ref{E:2})
and (\ref{E:4}) gives tree linear equations for $\nu_{2a}$,
$\nu_{2b}$ and  $\nu_{7a}$. Since $ \tiny{ \left|
\begin{matrix}
1&1&1\\
6&-2&1\\
13&1&0
\end{matrix}\right| } \not=0$, this system has an unique solution,
but this solution is not integral.

\noindent$\bullet$ Let $u$ be a unit of order $21$. By (\ref{E:1})
and Proposition \ref{P:4} we have  $\nu_{3a}+\nu_{7a}=1$. Using
Proposition \ref{P:1} for the characters $\chi_{2}, \chi_{3}$ of
$G$, we get the following system
\[
\begin{split}
\mu_{0}(u,\chi_{2},*) & = \textstyle \frac{1}{21} (48 \nu_{3a} + 12 \nu_{7a} + 36) \geq 0; \\ 
\mu_{7}(u,\chi_{2},*) & = \textstyle \frac{1}{21} (-24 \nu_{3a} - 6 \nu_{7a} + 24) \geq 0; \\ 
\mu_{1}(u,\chi_{3},*) & = \textstyle \frac{1}{21} (5 \nu_{3a} + 72) \geq 0, \\ 
\end{split}
\]
which has no nonnegative integral solution $(\nu_{3a},\nu_{7a})$.

\noindent$\bullet$ Let $u$ be a unit of order $22$. By (\ref{E:1})
and Proposition \ref{P:4} we have
$$
\nu_{2a}+\nu_{2b}+\nu_{11a}+\nu_{11b}=1.
$$

Put

\begin{equation}\label{E:9}
{ (\alpha,\beta)=\begin{cases}
(28,16), &\text{ if }\quad   \chi(u^{11}) = \chi(2a);\\
(20,24),\quad &\text{ if }\quad \chi(u^{11}) = \chi(2b);\\
(4,4),\quad &\text{ if }\quad  \chi(u^{11})=-2\chi(2a)+3\chi(2b);\\
(36,8),\quad &\text{ if }\quad  \chi(u^{11})=2\chi(2a)-\chi(2b);\\
(0,0),\quad &\text{ if }\quad  \chi(u^{11})=3\chi(2a)-2\chi(2b);\\
(12,32),\quad &\text{ if }\quad  \chi(u^{11})=-\chi(2a)+2\chi(2b).\\
\end{cases}}
\end{equation}
Moreover set
\begin{equation}\label{E:10}
t_1 = 3 \nu_{2a} -  \nu_{2b}, \qquad \ t_2 = 13 \nu_{2a}
+\nu_{2b},\qquad t_3 = 16 \nu_{2b} + 5 \nu_{11a} - 6 \nu_{11b}.
\end{equation}
Since  $|u^{11}|=2$ and $|u^{2}|=11$, for any character $\chi$ of
$G$ we need to consider $60$  cases, defined by parts (iii) and
(v) of the Theorem. We parameterize these six cases by values of
$(\alpha,\beta)$ from  (\ref{E:9}).

Case 1. Let $(\alpha,\beta)=(28,16)$. Then by Proposition
\ref{P:1} we get the system
\[
\begin{split}
\mu_{0}(u,\chi_{2},*) & = \textstyle \frac{1}{22} (20t_1 + \alpha)
\geq 0; \qquad
\mu_{11}(u,\chi_{2},*)  = \textstyle \frac{1}{22} (-20t_1 + \beta) \geq 0, \\ 
\end{split}
\]
that  has no integral solution.

Case 2. Let $(\alpha,\beta)=(20,24)$. Then by Proposition
\ref{P:1} we get the system
\begin{equation}\label{E:11}
\begin{split}
\mu_{0}(u,\chi_{2},*) & =\textstyle \frac{1}{22} (20t_1 +\alpha)
\geq 0; \quad
\mu_{11}(u,\chi_{2},*)  =\textstyle \frac{1}{22} (-20t_1 + \beta) \geq 0. \\ 
\end{split}
\end{equation}

Again by Proposition \ref{P:1} we obtain that
\begin{equation}\label{E:12}
\begin{split}
\mu_{0}(u,\chi_{3},*) & =\textstyle \frac{1}{22} (10 t_2 + 78)
\geq 0; \qquad
\mu_{1}(u,\chi_{3},*) =\textstyle \frac{1}{22} (t_2 + 76) \geq 0; \\ 
&\qquad \mu_{11}(u,\chi_{3},*)  =\textstyle \frac{1}{22} (-10 t_2 + 76) \geq 0; \\ 
&\mu_{1}(u,\chi_{14},*)  =\textstyle \frac{1}{22} (16 \nu_{2b} + 5 \nu_{11a} - 6 \nu_{11b} +\gamma ) \geq 0; \\ 
&\mu_{4}(u,\chi_{14},*)  =\textstyle \frac{1}{22} (- 16 \nu_{2b} - 5 \nu_{11a} + 6 \nu_{11b} + \delta) \geq 0, \\ 
\end{split}
\end{equation}
where
$$
\tiny{ (\gamma, \delta)= \begin{cases}
(886,918), &\text{ if  }\quad  \chi(u^{11}) = \chi(2b) \quad  \text{ and }\quad  \chi(u^{2})  = \chi(11a);\\
(875,907), &\text{ if  }\quad  \chi(u^{11}) = \chi(2b) \quad  \text{ and }\quad  \chi(u^{2}) = \chi(11b) ;\\
(930,962), &\text{ if  }\quad  \chi(u^{11}) = \chi(2b) \quad  \text{ and }\quad  \chi(u^{2}) =5\chi(11a)-4\chi(11b);\\
(853,885), &\text{ if  }\quad  \chi(u^{11}) = \chi(2b) \quad  \text{ and }\quad  \chi(u^{2}) = -2\chi(11a)+3\chi(11b);\\
(897,929), &\text{ if  }\quad  \chi(u^{11}) = \chi(2b) \quad  \text{ and }\quad  \chi(u^{2}) =  2\chi(11a)-\chi(11b);\\
(842,874), &\text{ if  }\quad  \chi(u^{11}) = \chi(2b) \quad  \text{ and }\quad  \chi(u^{2}) = -3\chi(11a)+4\chi(11b);\\
(831,863), &\text{ if  }\quad  \chi(u^{11}) = \chi(2a) \quad  \text{ and }\quad  \chi(u^{2}) = -4\chi(11a)+5\chi(11b);\\
(908,940), &\text{ if  }\quad  \chi(u^{11}) = \chi(2a) \quad  \text{ and }\quad  \chi(u^{2}) =  3\chi(11a)-2\chi(11b);\\
(864,896), &\text{ if  }\quad  \chi(u^{11}) = \chi(2a) \quad  \text{ and }\quad  \chi(u^{2}) = -\chi(11a)+2\chi(11b);\\
(919,951), &\text{ if  }\quad  \chi(u^{11}) = \chi(2b) \quad  \text{ and }\quad  \chi(u^{2}) = 4\chi(11a)-3\chi(11b).\\
\end{cases}}
$$
If we substitute the possible values of $(\gamma,\delta)$  into
(\ref{E:11})-- (\ref{E:12}), then it is easy to check that
$t_1=-1$ and $t_2=1$ (also we can calculate $t_3$). Now we
substitute back these values of $t_1$, $t_2$  and $t_3$ into
(\ref{E:10}). Then (\ref{E:1}) and (\ref{E:10}) gives four linear
equations for $\nu_{2a}$, $\nu_{2b}$, $\nu_{11a}$ and $\nu_{11b}$.
Since $ \tiny{ \left|
\begin{matrix}
1&1&1&1\\
3&-1&0&0\\
13&1&0&0\\
0&16&5&-6
\end{matrix}\right| } \not=0$, this system has an unique solution,
but this solution is not integral.

Case 3. Let $(\alpha,\beta)=(4,4)$. Then by Proposition \ref{P:1}
we get the system
\[
\begin{split}
\mu_{0}(u,\chi_{2},*) & =\textstyle \frac{1}{22} (20t_1 + 4) \geq
0; \qquad
\mu_{2}(u,\chi_{2},*)  =\textstyle \frac{1}{22} (-2t_1 + 4) \geq 0; \\ 
\mu_{0}(u,\chi_{3},*) & =\textstyle \frac{1}{22} (10t_2 + 54) \geq
0; \quad
\mu_{11}(u,\chi_{3},*) =\textstyle \frac{1}{22} (- 10 t_2 + 100) \geq 0; \\ 
\mu_{1}(u,\chi_{14},*) & =\textstyle \frac{1}{22} (t_3 + \gamma)
\geq 0; \qquad\;
\mu_{4}(u,\chi_{14},*)  =\textstyle \frac{1}{22} (- t_3 + \delta) \geq 0, \\ 
\end{split}
\]
where
$$
\tiny{ (\gamma, \delta)= \begin{cases}

(854,950), &\text{ if }\quad    \chi(u^{2}) = \chi(11a);\\
(843,939), &\text{ if  }\quad   \chi(u^{2}) = \chi(11b);\\
(898,994), &\text{ if  }\quad   \chi(u^{2}) = 5\chi(11a)-4\chi(11b);\\
(821,917), &\text{ if  }\quad  \chi(u^{2}) = -2\chi(11a)+3\chi(11b);\\
(865,961), &\text{ if  }\quad   \chi(u^{2}) = 2\chi(11a)-\chi(11b);\\
(810,906), &\text{ if  }\quad   \chi(u^{2}) = -3\chi(11a)+4\chi(11b);\\
(799,895), &\text{ if  }\quad  \chi(u^{2}) = -4\chi(11a)+5\chi(11b);\\
(876,972), &\text{ if  }\quad  \chi(u^{2}) = 3\chi(11a)-2\chi(11b);\\
(832,928), &\text{ if  }\quad   \chi(u^{2}) = -\chi(11a)+2\chi(11b);\\
(887,983), &\text{ if  }\quad   \chi(u^{2}) = 4\chi(11a)-3\chi(11b).\\
\end{cases}
}
$$
By easy calculation we obtain that $t_1=2$ and $t_2\in\{-1,10\}$.
This case is similar to the previous ones, so we can conclude that
there is no integral solution.

Case 4. Let $(\alpha,\beta)=(36,8)$. Then by Proposition \ref{P:1}
we get the system
\[
\begin{split}
\mu_{0}(u,\chi_{2},*) & =\textstyle \frac{1}{22} (20t_1 +\alpha)
\geq 0; \qquad
\mu_{11}(u,\chi_{2},*)  =\textstyle \frac{1}{22} (-20t_1 +\beta) \geq 0, \\ 
\end{split}
\]
which leads to a contradiction.

Case 5. Let $(\alpha,\beta)=(0,0)$. Then by Proposition \ref{P:1}
we get the system
\[
\begin{split}
\mu_{1}(u,\chi_{2},*) & =\textstyle \frac{1}{22} (2t_1 ) \geq 0;
\qquad\qquad\quad
\mu_{11}(u,\chi_{2},*)  =\textstyle \frac{1}{22} (-20t_1 ) \geq 0; \\ 
\mu_{0}(u,\chi_{3},*) & =\textstyle \frac{1}{22} (10t_2 + 114)
\geq 0; \qquad
\mu_{11}(u,\chi_{3},*)  =\textstyle \frac{1}{22} (-10 t_2 + 40) \geq 0; \\ 
\mu_{1}(u,\chi_{14},*) & =\textstyle \frac{1}{22} (t_3 + \gamma)
\geq 0; \qquad\qquad
\mu_{4}(u,\chi_{14},*)  =\textstyle \frac{1}{22} (- t_3 +\delta) \geq 0, \\ 
\end{split}
\]
where
$$
\tiny{ (\gamma, \delta)= \begin{cases} (935,880) &\text{ if}\quad \chi(u^{2}) = \chi(11a) ;\\
(923,859), &\text{ if  }\quad \chi(u^{2}) = \chi(11b);\\
(978,914), &\text{ if  }\quad\chi(u^{2}) = 5\chi(11a)-4\chi(11b);\\
(901,837), &\text{ if  }\quad\chi(u^{2}) = -2\chi(11a)+3\chi(11b);\\
(945,881), &\text{ if  }\quad\chi(u^{2}) = 2\chi(11a)-\chi(11b);\\
(890,826), &\text{ if  }\quad\chi(u^{2}) = -3\chi(11a)+4\chi(11b);\\
(879,815), &\text{ if  }\quad \chi(u^{2}) = -4\chi(11a)+5\chi(11b);\\
(956,892), &\text{ if  }\quad\chi(u^{2}) = 3\chi(11a)-2\chi(11b);\\
(912,848), &\text{ if  }\quad\chi(u^{2}) = -\chi(11a)+2\chi(11b);\\
(967,903), &\text{ if  }\quad\chi(u^{2}) = 4\chi(11a)-3\chi(11b).\\
\end{cases}}
$$
By easy calculation we obtain that $t_1=0$ and $t_2\in\{-7,4\}$.
This case is similar to the Case 3, so we can conclude that there
is no integral solution in this case too.

Case 6. Finally, let $(\alpha,\beta)=(36,8)$. By Proposition
\ref{P:1} we obtain  the system:
\[
\begin{split}
\mu_{0}(u,\chi_{2},*) & =\textstyle \frac{1}{22} (20t_1 + 12) \geq
0; \qquad
\mu_{11}(u,\chi_{2},*) =\textstyle \frac{1}{22} (-20t_1 + 32) \geq 0, \\ 
\end{split}
\]
which has no integral solution.

\noindent$\bullet$ Let $u$ be a unit of order $33$. By (\ref{E:1})
and Proposition \ref{P:4} we get $ \nu_{3a}+\nu_{11a}+\nu_{11b}=1.
$ Again, using Proposition \ref{P:1} we obtain that
\[
\begin{split}
\mu_{0}(u,\chi_{2},*) & = \textstyle \frac{1}{33} (80 \nu_{3a} +
30) \geq 0; \quad
\mu_{11}(u,\chi_{2},*)  = \textstyle \frac{1}{33} (-40 \nu_{3a} + 18) \geq 0, \\ 
\end{split}
\]
that has no integral solution.

\noindent$\bullet$ Let $u$ be a unit of order $35$. By (\ref{E:1})
and Proposition \ref{P:4} we have that
$$
    \nu_{5a}+\nu_{5b}+\nu_{5c}+\nu_{7a}=1.
$$
Since  $|u^7|=5$, for any character $\chi$  we need to consider
$23$ cases, defined by part (iv) of the Theorem. Using Proposition
\ref{P:1}, we divide these $23$ cases into five groups:

Group 1. Let\quad  $\chi(u^{7})$ belongs to the following set
\[
\begin{split}
\{\; & \chi(5a)-\chi(5b)+\chi(5c),\quad \chi(5a)+4\chi(5b)-4\chi(5c),\\
&\chi(5a)+3\chi(5b)-3\chi(5c), \quad \chi(5a)-2\chi(5b)+2\chi(5c),\\
& \chi(5a)+2\chi(5b)-2\chi(5c),\quad
\chi(5a)+\chi(5b)-\chi(5c),\quad \chi(5a)\;\}.
\end{split}
\]
Applying Proposition \ref{P:1} to the character $\chi_{2}$ we
construct the following system
\[
\begin{split}
\mu_{5}(u,\chi_{2},*) & = \textstyle \frac{1}{35} (12 \nu_{5a} - 8 \nu_{5b} - 8 \nu_{5c} - 4 \nu_{7a} + 9) \geq 0; \\ 
\mu_{0}(u,\chi_{2},*) & = \textstyle \frac{1}{35} (-6(12 \nu_{5a} - 8 \nu_{5b} - 8 \nu_{5c} - 4 \nu_{7a}) + 16) \geq 0, \\ 
\end{split}
\]
which has no integral solution.

\noindent Group  2. Let\quad $\chi(u^{7})$ belongs to the
following set
\[
\begin{split}
\{\; \chi(5b),\quad \chi(5c),&\quad -\chi(5b)+2\chi(5c),\quad
3\chi(5b)-2\chi(5c),\\
&\qquad 2\chi(5b)-\chi(5c),\quad -2\chi(5b)+3\chi(5c)\;\}.
\end{split}
\]
Using Proposition \ref{P:1} to the character $\chi_{2}$ we get the
system
\[
\begin{split}
\mu_{7}(u,\chi_{2},*) & = \textstyle \frac{1}{35} (18 \nu_{5a} - 12 \nu_{5b} - 12 \nu_{5c} - 6 \nu_{7a} + 26) \geq 0; \\ 
\mu_{0}(u,\chi_{2},*) & = \textstyle \frac{1}{35} (-72 \nu_{5a} + 48 \nu_{5b} + 48 \nu_{5c} + 24 \nu_{7a} + 36) \geq 0, \\ 
\end{split}
\]
which has no integral solution.

\noindent Group  3. Let\quad  $\chi(u^{7})$ belongs to the
following set
\[
\begin{split}
\{ & \quad -2\chi(5a)-\chi(5b)+4\chi(5c),\quad
-2\chi(5b)+3\chi(5c),\\
&\quad  -2\chi(5a)-2\chi(5b)+5\chi(5c),\quad
-2\chi(5a)+\chi(5b)+2\chi(5c)\quad \}.
\end{split}
\]
By Proposition \ref{P:1} we have the following system of
inequalities
\[
\begin{split}
\mu_{7}(u,\chi_{2},*) & = \textstyle \frac{1}{35} (18 \nu_{5a} - 12 \nu_{5b} - 12 \nu_{5c} - 6 \nu_{7a} + 16) \geq 0; \\ 
\mu_{0}(u,\chi_{2},*) & = \textstyle \frac{1}{35} (-72 \nu_{5a} + 48 \nu_{5b} + 48 \nu_{5c} + 24 \nu_{7a} + 76) \geq 0, \\ 
\end{split}
\]
which has no integral solution.

\noindent Group  4. Let\quad  $\chi(u^{7})$ belongs to the
following set
\[
\begin{split}
\{  \; -\chi(5a)-&\chi(5b)+3\chi(5c),\quad
-\chi(5a)+2\chi(5c),\quad
-\chi(5a)+2\chi(5b),\\
&\quad -\chi(5a)-2\chi(5b)+4\chi(5c),\quad
-\chi(5a)+\chi(5b)+\chi(5c)\; \}.
\end{split}
\]
Using Proposition \ref{P:1} to the character $\chi_{2}$ we get the
system
\[
\begin{split}
\mu_{7}(u,\chi_{2},*) & = \textstyle \frac{1}{35} (18 \nu_{5a} - 12 \nu_{5b} - 12 \nu_{5c} - 6 \nu_{7a} + 21) \geq 0; \\ 
\mu_{0}(u,\chi_{2},*) & = \textstyle \frac{1}{35} (-72 \nu_{5a} + 48 \nu_{5b} + 48 \nu_{5c} + 24 \nu_{7a} + 56) \geq 0, \\ 
\end{split}
\]
which has no integral solution.


\noindent Group 5. Finally, let \quad
$\chi(u^{7})=-3\chi(5a)+4\chi(5c)$. By Proposition \ref{P:1} we
get
\[
\begin{split}
\mu_{0}(u,\chi_{3},*) & = \textstyle \frac{1}{35} (48 \nu_{5a} - 72 \nu_{5b} + 48 \nu_{5c} + 85) \geq 0; \\ 
\mu_{7}(u,\chi_{3},*) & = \textstyle \frac{1}{35} (-12 \nu_{5a} + 18 \nu_{5b} - 12 \nu_{5c} + 75) \geq 0,  \\ 
\end{split}
\]
which has no integral solution too.


\noindent$\bullet$ Let $u$ be a unit of order $55$. By (\ref{E:1})
and Proposition \ref{P:4} we have
$$
\nu_{5a}+\nu_{5b}+\nu_{5c}+\nu_{11a}+\nu_{11b}=1.
$$
Since  $|u^{11}|=5$ and $|u^{5}|=11$, for any character $\chi$ of
$G$ we need to consider $230$ cases, defined by parts (iv) and (v)
of the Theorem. Using our implementation of the Luthar--Passi method,
which we intended to make available in the GAP package LAGUNA \cite{LAGUNA},
we can employ Proposition \ref{P:1} to construct in all 230 cases the 
systems of inequalities. Actually in all cases we obtain a system of two
inequalities, and a lot of cases lead to the same system.
We present here two cases that yield the same system.

Let \; $\chi(u^{11}) = \chi(5a)$ \; and either \; $\chi(u^{5}) =
\chi(11a)$\; or\; $\chi(u^{5}) = 5\chi(11a)-4\chi(11b)$. By
Proposition \ref{P:1} we obtain the following system of
inequalities
\[
\begin{split}
\mu_{5}(u,\chi_{2},*) & = \textstyle \frac{1}{55} (12 \nu_{5a} - 8 \nu_{5b} - 8 \nu_{5c} + 10) \geq 0; \\ 
\mu_{0}(u,\chi_{2},*) & = \textstyle \frac{1}{55} (-10(12\nu_{5a} - 8 \nu_{5b} -8 \nu_{5c}) + 10) \geq 0, \\ 
\end{split}
\]
which has no integral solution  such that all
$\mu_{i}(u,\chi_{j},p)$ are nonnegative integers.

\noindent$\bullet$ Let $u$ be a unit of order $77$. By (\ref{E:1})
and Proposition \ref{P:4} we obtain  that
$$
\nu_{7a}+\nu_{11a}+\nu_{11b}=1.
$$
Finally, using Proposition \ref{P:1} we get the system of
inequalities:
\[
\begin{split}
\mu_{0}(u,\chi_{2},*) & = \textstyle \frac{1}{77} (60 \nu_{7a} + 28) \geq 0; \\ 
\mu_{0}(u,\chi_{3},5) & = \textstyle \frac{1}{77} (-60 \nu_{7a} + 49) \geq 0, \\ 
\end{split}
\]
which has no integral solution such that all
$\mu_{i}(u,\chi_{j},p)$ are nonnegative integers.

\bibliographystyle{plain}
\bibliography{HS_Kimmerle}

\end{document}